\documentclass[a4paper,11pt]{article}
\usepackage{amsmath}
\usepackage{amssymb}
\usepackage{mathrsfs}
\usepackage{amsfonts}
\usepackage{amsthm}
\usepackage{pstricks, pst-node, pst-text, pst-3d,psfrag}
\usepackage{graphicx}
\usepackage{indentfirst}
\usepackage{enumerate}
\usepackage{subfigure}
\usepackage{cite}
\usepackage[colorlinks=true]{hyperref}
\hypersetup{urlcolor=blue, citecolor=red}

\theoremstyle{plain}
\newtheorem{thm}{Theorem}[section]
\newtheorem{thmm}{Theorem}[subsection]
\newtheorem*{conjecture}{Conjecture}
\newtheorem*{reconjecture}{Conjecture A}

\newtheorem{prop}[thm]{Proposition}

\newtheorem{lem}[thm]{Lemma}
\theoremstyle{definition}
\newtheorem{defn}[thm]{Definition}

\newtheorem{rmk}[thm]{Remark}

\numberwithin{equation}{section}

\makeatletter 
\@addtoreset{equation}{section}
\makeatother  

\newcommand\abs[1]{\lvert#1\rvert}

\begin{document}

\title{Almost sure convergence of differentially positive systems on a globally orderable manifold}

\setlength{\baselineskip}{16pt}

\author {
Lin Niu\thanks{Supported by the National Natural Science Foundation of China No.12201034 and 12090012.}\\
School of Mathematics and Physics\\
University of Science and Technology Beijing\\
Beijing, 100083, P. R. China
\\[3mm]
Yi Wang\thanks{Supported by the National Key R\&D Program of China (2024YFA1013603, 2024YFA1013600), NSFC No.12331006, and Strategic Priority Research Program of Chinese Academy of Sciences (XDB0900100).}
\\
School of Mathematical Sciences\\
University of Science and Technology of China\\
Hefei, Anhui, 230026, P. R. China
\\[3mm]
Yufeng Zhang\thanks{Supported by the National Natural Science Foundation of China No.12401240.}
\\
Department of Mathematical and Statistical Sciences\\ University of Alberta\\ Edmonton, AB T6G 2G1, Canada
}

\date{}

    \maketitle

\begin{abstract}
	Differentially positive systems are nonlinear systems whose linearization along trajectories preserves a cone field on a smooth Riemannian manifold. The structures of cone field come from general relativity and Lie theory. We prove that on a globally orderable manifold, the set of convergent points has full Riemann–Lebesgue measure, thus establishing almost sure convergence. This result thereby resolves a measure-theoretic form of the conjecture posed by Forni and Sepulchre in 2016 for such manifolds.

\vskip 3mm

\end{abstract}


\section{Introduction}

The present paper investigates a class of nonlinear dynamical systems whose linearizations along trajectories preserve a cone field. Such systems are referred to as differentially positive systems. Differential analysis offers a general framework for studying nonlinear dynamical systems by examining their linearizations at each point in the state space, motivated by the observation that local behavior can strongly influence global dynamics. Positivity refers to the property of preserving a cone in a vector space. The differential perspective extends this notion to the nonlinear setting by introducing a cone field on a manifold $M$. Specifically, a cone field $C_M$ assigns to each point $x\in M$ a closed convex cone $C_M(x)$ in the tangent space $T_{x}M$.

The concept of differential positivity was first introduced by Forni and Sepulchre \cite{Forni2015,ForniandSepulchre14,ForniandSepulchre16} in their study of the dynamics of the pendulum model. In \cite[Section VIII]{ForniandSepulchre16}, the damped pendulum is written as
\begin{equation*}
	\begin{cases}
		\frac{dx}{dt} = y,\\
		\frac{dy}{dt} =-\sin x-ky+u,
	\end{cases}
\end{equation*}
where $(x, y)\in M= \mathbb{S}^1\times \mathbb{R}^1$, $u$ is the torque input and $k\geq 0$ is the damping coefficient. The pendulum is differentially positive for $k\geq 2$ with respect to the cone field $C_M(x,y)=\{(v_1, v_2)\in T_{(x,y)}M \colon v_1\geq 0, v_1+v_2\geq 0\}$.
Meanwhile, the application of differential positivity to a class of consensus protocols on the $N$-torus $\mathbb{T}^N$ is presented in \cite[Section 4]{MostajeranSepulchre18SIAM}, where the system $$\frac{d\theta_k}{dt}=\sum\limits_{i:(k,i)\in \mathcal{E}}g_{ki}(\theta_i-\theta_k)$$ on $\mathbb{T}^N$ for a strongly connected communication graph $(\mathcal{V},\mathcal{E})$ is differentially positive with respect to the cone field $C_{\mathbb{T}^N}(\theta)=\{\delta \theta\in T_{\theta}\mathbb{T}^N \colon \delta \theta^i\geq 0, \delta \theta=\sum\limits_{i}\delta \theta^i \frac{\partial}{\partial \theta^i}\}$ on the set $\mathbb{T}^N_{\pi}=\{\theta\in \mathbb{T}^N \colon \abs{\theta_k-\theta_i}<\pi, (i,k)\in \mathcal{E}\}$, where $\theta_k\in \mathbb{S}^1$, $\theta=(\theta_1,\cdots,\theta_N)$. Here, the $N$-tuple of vector fields $(\frac{\partial}{\partial \theta^1},\cdots,\frac{\partial}{\partial \theta^N})$  are the standard smooth global frame on $\mathbb{T}^N$, and each function $g_{ki}$ is an odd, $2\pi$-periodic, continuously differentiable on $(-\pi,\pi)$ satisfying $g_{ki}(0)=0$, $g'_{ki}(\alpha)>0$ for any $\alpha\in (-\pi,\pi)$.
This consensus protocol generalizes the Kuramoto model (see \cite[Section V]{Acebronetc2005}), a widely used framework for analyzing synchronization phenomena in smart grids and engineered complex oscillator networks (see, e.g., \cite{Dorfleretc2013,Simpsonetc13}).
One may refer to more examples of differentially positive systems in \cite{Forni2015,ForniandSepulchre14,ForniandSepulchre16,MostajeranSepulchre18SIAM,
MostajeranSepulchre18homogeneous}.

Forni and Sepulchre \cite{ForniandSepulchre16} first investigated the dynamics of differentially positive systems on a Riemannian manifold $M$. Under the assumption that every orbit of the system is defined for all $t \in \mathbb{R}$ with compact closure, they established a dichotomy which characterizes the structure of the limit sets of differentially positive systems (see \cite[Theorem 4]{ForniandSepulchre16}). Based on this, Forni and Sepulchre \cite[p.353]{ForniandSepulchre16} further posed the following

\begin{conjecture}[Forni and Sepulchre \cite{ForniandSepulchre16}]
	 For almost every $x\in M$, the $\omega$-limit set $\omega(x)$ is given by either a fixed point, or a limit cycle, or fixed points and connecting arcs.
\end{conjecture}

This conjecture concerns the asymptotic behavior of orbits as time tends to $+\infty$, aiming to characterize what typically occurs for ``most" initial conditions (see more details in \cite[p.353]{ForniandSepulchre16}).

In our previous work \cite{NW24}, two of the present authors tackled this conjecture and have established the asymptotic behavior of generic (i.e., ``almost all" in the topological sense) initial points on a globally orderable manifold. Nevertheless, there are many important situations in which properties that hold ``generically" (in the topological sense) differ from those that hold ``almost surely" (in the measure-theoretic sense), and vise-versa. As a matter of fact, these two ``typical" notions are indeed quite distinct (see, e.g., \cite{HuntandKaloshin10,WangandYaoandZhangJDDE22} and references therein).

\vskip 2mm

In the present paper, we will focus on this conjecture {\it in the measure-theoretic sense}, and attempt to establish the asymptotic behavior of ``almost every" orbit for differentially positive systems.

\vskip 2mm

For this purpose, we begin by recalling some insight for the cone fields on a smooth manifold. One of the earliest examples of cone fields arises in general relativity, where a time-orientable space-time naturally determines, at each point, a Lorentzian cone in the tangent space, i.e., a closed, convex, pointed cone representing all future-directed non-spacelike directions (see, e.g., \cite{Beemandehrlich81, Hawkingandellis73, Penrose72}). Such cone field defines the class of non-spacelike (or causal) curves, which are precisely those curves whose tangent vectors lie in the cone at each point along the curve.
In this setting, each point of space-time corresponds to an event, and a signal can be sent from $p$ to $q$ if there exists a future-directed causal curve from $p$ to $q$. That is, causally related points in space-time are connected by a non-spacelike curve determined by the Lorentzian cone field.

From this perspective, a cone field on the manifold $M$ naturally induces a ``conal order relation" as follows: two points $x,y\in M$ are {\it conal ordered}, denoted by $x\leq_M y$, if there exists a conal curve on $M$ starting at $x$ and ending at $y$. Here, a conal curve is a continuous and piecewise continuously differentiable curve whose tangent vector lies in the cone at every point along the curve wherever the derivative is defined (see Definition \ref{conal curve}). In particular, in the context of space-time geometry (see \cite{Beemandehrlich81, Hawkingandellis73, Penrose72}), causal curves correspond exactly to conal curves, and causally related points are precisely those that are conal ordered. Moreover, conal order has notable applications beyond relativity, including in the theory of hyperbolic partial differential equations and harmonic analysis (see \cite{Faraut87, Faraut91}), as well as in the study of Wiener-Hopf operators on symmetric spaces \cite{HilgertNeeb}.

As one may know, unlike the standard partial order induced by a closed convex cone in a vector space, the conal order relation ``$\leq_M$" induced by a cone field on $M$ is reflexive and transitive, but \textit {not necessarily antisymmetric}. As a matter of fact, there exist examples of manifolds $M$ that admit closed conal curves, such as closed timelike curves in time-orientable space-times (see \cite[Chapter 5]{Hawkingandellis73}), which demonstrate the failure of antisymmetry. More importantly, the relation ``$\leq_M$" is \textit {not necessarily closed}; that is, the set $\{(x,y)\in M\times M: x\leq_M y \}$ need not be closed in $M\times M$ (see Lawson \cite[p.299]{Lawson89}, or Neeb \cite[p.470]{Neeb91}). For instance, such failure occurs in Minkowski space (see Hawking and Ellis \cite[p.183]{Hawkingandellis73} or Penrose \cite[p.12]{Penrose72}).

\vskip 2mm

From the viewpoint of order relations, the flow of a differentially positive system naturally preserves the conal order ``$\leq_M$" (see Proposition \ref{strong order preserving}). As a first attempt to solve the measure-theoretic version of the conjecture, we introduce the following assumption on the manifold $M$:

\begin{enumerate}[{\bf (H1)}]
	\item $M$ is globally orderable equipped with a continuous solid cone field.
\end{enumerate}

(H1) means that the conal order ``$\le_M$" on $M$ is a partial order (i.e., it is additionally antisymmetric, see Definition \ref{globally orderable}). The notion of globality for the conal order ``$\leq_M$" has its roots in the Lie theory of semigroups, where it plays a central role in the study of invariant cone structures on homogeneous spaces (see \cite{Lawson89, HilgertHofmannLawson, HilgertNeeb93, Neeb91, NW24}).
A fundamental question explored in \cite{Lawson89} is whether the conal order induced by a cone field can be extended to a genuine partial order on the entire manifold. Lawson \cite[Section 5]{Lawson89} established an equivalence between the globality of the conal order on a homogeneous manifold and the globality of the associated Lie wedge in the Lie group acting on it (see also in Neeb \cite[Theorem 1.6]{Neeb91}, or \cite[Section 4.3]{HilgertNeeb93}).

In fact, (H1) arises naturally in various settings. A prominent example is the homogeneous space of positive definite matrices equipped with the affine-invariant cone field (see \cite[Section 3]{MostajeranSepulchre18orderingmatrices}). Moreover, (H1) precludes the existence of closed conal curves (see \cite[Section 5]{Lawson89}). In the context of space-times, such closed conal (or causal) curves are known to lead to paradoxes related to causality and are commonly said to ``violate causality" (\cite{Beemandehrlich81, Hawkingandellis73, Hawkingandsachs74, Penrose72}).

\vskip 2mm

As a consequence, by (H1), the occurrence of closed conal curves, such as limit cycles, fixed points and connecting arcs, can be excluded. In this way, the Forni-Sepulchre Conjecture naturally reduces to the following formulation

\begin{reconjecture}
	If {\rm (H1)} holds, then for almost every $x\in M$, the $\omega$-limit set $\omega(x)$ is a singleton.
\end{reconjecture}


In the present paper, we will prove Conjecture A in the measure-theoretic sense, under the following two reasonable assumptions:

\begin{enumerate}
	\item[{\bf (H2)}] The conal order ``$\leq_M$" is quasi-closed.
	\item[{\bf (H3)}] Both the cone field $C_M$ and the Riemannian metric on $M$ are $\Gamma$-invariant.
\end{enumerate}

The quasi-closed order relationship ``$\leq_M$" in (H2) is inspired by the causal continuity of space-times in Hawking and Sachs \cite{Hawkingandsachs74} (see also \cite{NW24} for more details). It means that ``$x\leq_M y$ whenever $x_n\to x$ and $y_n\to y$ as $n\to \infty$ and $x_n\ll_M y_n$ for all integer $n\ge 0$" (see Definition \ref{quasi-closed}). Here, we write $x\ll_M y$ if there exists a so-called \textit{strictly} conal curve (whose tangent vector lies in the interior of cone at every point along the curve) on $M$ starting at $x$ and ending at $y$ (see Definition \ref{conal curve}).

As for (H3), a cone field $C_M$ is called \textit{$\Gamma$-invariant} if there exists a linear invertible mapping $\Gamma(x_1,x_2) : T_{x_1}M \to T_{x_2}M$ for all $x_1,x_2\in M$, such that $\Gamma(x_1,x_2)C_M(x_1)=C_M(x_2)$.
While, the Riemannian metric $(\cdot,\cdot)_x$ on $M$ is \textit{$\Gamma$-invariant} if $(u,v)_{x_1}=(\Gamma(x_1,x_2)u,\Gamma(x_1,x_2)v)_{x_2}$ for all $x_1,x_2\in M$ and $u,v\in T_{x_1}M$.
Actually, $\Gamma$-invariance is motivated by the homogeneous structure of the manifolds (\cite{NW24}).

Our main result is the following


\renewcommand{\thethmm}{\Alph{thmm}}

\begin{thmm}[Almost Sure Convergence]\label{genericconvergence}
	Assume that {\rm (H1)-(H3)} hold. Then the set of convergent points has full Riemann-Lebesgue measure.
\end{thmm}

Theorem \ref{genericconvergence} will be proved in Section $4$ (see Theorem \ref{measure zero}). It concludes the fact that ``the typical orbit will converge to an equilibrium" holds almost-surely in the measure-theoretic sense.

In our present work, $M$ is a Riemannian manifold which gives rise to a canonical measure structure (namely, the Riemann-Lebesgue measure $\lambda_M$). A subset $\mathcal{C}\subset M$ has full Riemann–Lebesgue measure if the complement of $\mathcal{C}$ in $M$ is a $\lambda_M$-null set. In fact, there are various equivalent ways of describing what it means for a subset in $M$ to be $\lambda_M$-null (see, e.g., Proposition 1.6 in \cite[Chapter XII]{AmannandEscher}). What one needs here is a description of the $\lambda_M$-null set in terms of local charts. We select one such description (i.e., the concept of sets of measure zero, in Definition \ref{measurezerodefi} and Remark \ref{measueRMK}) that will be convenient for our purposes in this paper. Here, a point $x$ is a convergent point if the $\omega$-limit set $\omega(x)$ is a singleton.


It deserves to point out that our Theorem \ref{genericconvergence} is built upon the following {\it distinctive phenomenon}: For any simply ordered set (i.e., any two distinct points on it are ordered with respect to ``$\le_M$"), at most countably many points on it are non-convergent (see Lemma \ref{countable set}).
In fact, in our previous work \cite{NW24}, a key insight for the ``generic convergence" is to show that if $x\notin \mathcal{C}$ (here, the set of convergent points is denoted by $\mathcal{C}$), then it must belong to the closure of $\mathcal{C}$ (i.e., $x\in \overline{\mathcal{C}}$, see \cite[Theorem 4.1]{NW24}). In other words, any such $x$ can be approached topologically by the points in the interior of $\mathcal{C}$. However, we stress that this key-point is not enough to guarantee the ``almost sure convergence" at all. In order to prove our main Theorem \ref{genericconvergence} here, we need the aformentioned phenomenon (Clearly, such phenomenon implies $x\in \overline{\mathcal{C}}$ automatically). Together with Fubini's Theorem, we succeed in establishing that the set of non-convergent points in each local chart has measure zero (see Remark \ref{rkofthm} for more details). Thus, the set of non-convergent points has measure zero in $M$.

\vskip 2mm

Differential positivity, as introduced by Forni and Sepulchre \cite{ForniandSepulchre16}, reduces to classical monotonicity theory when considered with respect to a constant cone field on a flat space $M$, i.e., a closed convex cone in a vector space (see \cite{H88,H85,HS05,S95,S17,Polacik02,HessandPolacik,SmithandThieme91,WangandYaoandZhangPAMS22,WangandYaoandZhangJDDE22}). From this perspective, differentially positive systems can be viewed as a natural extension of classical monotone dynamical systems to nonlinear manifolds. Consequently, Theorem \ref{genericconvergence} recovers the classical Hirsch's prevalent convergence Theorem (\cite{EncisoHirschSmith08,WangandYaoandZhangJDDE22}).

The paper is organized as follows. In Section $2$, we introduce some notations and standard assumptions. And we summarize some preliminary results in Section $3$. Theorem \ref{genericconvergence} (i.e., Theorem \ref{measure zero}) with its proof will be given in Section $4$. 

\section{Notations and standard hypotheses}

Throughout this paper, let $M$ denote a smooth manifold of dimension $r$. The tangent bundle is denoted by $TM$ and the tangent space at a point $x\in M$ by $T_{x}M$. The manifold $M$ is equipped with a Riemannian metric, represented by an inner product $(\cdot,\cdot)_x$ on each tangent space $T_{x}M$. For any $v\in T_{x}M$, we define $\abs{v}_x:=\sqrt{(v,v)_x}$, and, when no confusion arises, we omit the subscripts $x$. The Riemannian metric induces a distance function $d$ on $M$, and we assume that $(M,d)$ is a complete metric space. 


Let $E$ be a $r$ dimensional real linear space. A non-empty closed subset $K$ of the linear space $E$ is called a \textit{closed convex cone} if (i) if $x\in K$ and $\alpha \geq 0$, then $\alpha x\in K$; (ii) if $x,y\in K$, then $x+y\in K$; (iii) $K \cap (-K) = \{0\}$. A convex cone is \textit{solid} if the interior $\text{Int} K$ of $K$ is non-empty. A convex cone $K'$ in $E$ is said to \textit{surround} a cone $K$ if $K\setminus \{ 0\}$ is contained in the interior of $K'$ (see, e.g., \cite{Lawson89}).

\begin{defn}
	A \textit{cone field} is a map $x \mapsto C_{M}(x)$ on a manifold $M$, such that $C_{M}(x)$ is a convex cone in $T_{x}M$ for each $x\in M$.
\end{defn}
A manifold endowed with a cone field is called a \textit{conal manifold}. The cone field is said to be solid if each convex cone $C_{M}(x)$ is solid.

A cone field $C_M$ is called \textit{$\Gamma$-invariant} if there exists a linear invertible map $\Gamma(x_1,x_2) : T_{x_1}M \to T_{x_2}M$ for each $x_1,x_2\in M$, such that $\Gamma(x_1,x_2)C_M(x_1)=C_M(x_2)$. Moreover, $\Gamma$ is continuous in $(x_1, x_2)$ and satisfies $\Gamma(x,x)=\text{Id}_x$ for all $x\in M$, where $\text{Id}_x : T_xM \to T_xM$ denotes the identity map.
The Riemannian metric is said to be \textit{$\Gamma$-invariant} if $(u,v)_{x_1}=(\Gamma(x_1,x_2)u,\Gamma(x_1,x_2)v)_{x_2}$ for all $x_1,x_2\in M$ and $u,v\in T_{x_1}M$. In this case, $\abs{v}_{x_1}=\abs{\Gamma(x_1,x_2) v}_{x_2}$ holds for all $x_1,x_2\in M$ and $v\in T_{x_1}M$.

\begin{rmk}
	The notion of $\Gamma$-invariance is motivated by the homogeneous structure of manifolds (see \cite{NW24}). A classical example is a homogeneous space $M=G/H$ of a connected Lie group $G$, where each point $x\in M$ is assigned a closed convex cone in the tangent space $T_{x}M$ such that the cone field is invariant under the action of $G$ (see, e.g., \cite{Lawson89,Neeb91,HilgertNeeb93}); and moreover, any homogeneous Riemannian metric on such a space is naturally $\Gamma$-invariant.
\end{rmk}

Let $\Phi : U \to V$ be a smooth chart, where $U\subset M$ is an open set containing $x$ and $V\subset \mathbb{R}^{r}$ is open. For any $y\in U$, denote by $C_{M}^{\Phi}(y)\subset \mathbb{R}^r$ the representation of the cone field in this chart, i.e., $d\Phi(y)(C_{M}(y))=C_{M}^{\Phi}(y)$.
A cone field $x \mapsto C_{M}(x)$ on $M$ is said to be \textit{upper semicontinuous at} $x\in M$ if given a smooth chart $\Phi : U \to \mathbb{R}^{r}$ at $x$ and a convex cone $C'$ surrounding $C_{M}^{\Phi}(x)$, there exists a neighborhood $W$ of $x$ such that $C_{M}^{\Phi}(y)\subset C'$ for all $y\in W$. The cone field is said to be \textit{lower semicontinuous at} $x$ if given any open set $N$ such that $N\cap C_{M}^{\Phi}(x)$ is non-empty, there exists a neighborhood $W$ of $x$ such that $N\cap C_{M}^{\Phi}(y)$ is non-empty for all $y\in W$. A cone field is \textit{continuous at $x$} if it is both lower and upper semicontinuous at $x$ and \textit{continuous} if it is continuous at every $x$.

\begin{defn}\label{conal curve}
    A continuous and piecewise continuously differentiable curve $t \mapsto \gamma(t)$ defined on $[t_0,t_1]$ into a conal manifold $M$ (equipped with a continuous cone field $C_M$) is a \textit{conal curve} if the derivative $\gamma'(t)$ is in $C_{M}(\gamma(t))$ whenever $t_0\leq t<t_1$, in which the derivative is the right hand derivative at those finitely many points where the derivative is not continuous. Moreover, $\gamma$ is a \textit{strictly conal curve} if $\gamma$ is a conal curve and the derivative $\gamma'(t)$ is in $\text{Int} C_{M}(\gamma(t))$ for $t_0\leq t<t_1$.
\end{defn}

We say that two points $x,y\in M$ are ordered, denoted by $x\leq_M y$, if there is a conal curve $\alpha \colon [t_0,t_1]\subset \mathbb{R} \to M$ with $\alpha(t_0)=x$ and $\alpha(t_1)=y$.
This relation defines an \textit{order} on $M$, which is always reflexive (i.e., $x\leq_M x$ for all $x\in M$) and transitive (i.e., $x\leq_M y$ and $y\leq_M z$ implies $x\leq_M z$). The order ``$\leq_M$" is referred to as the \textit{conal order}. We write $x<_M y$ if $x\leq_M y$ and $x\neq y$.
We also write $x\ll_M y$ if there is a strictly conal curve $\gamma$ with $\gamma(t_0)=x$ and $\gamma(t_1)=y$. Clearly, the relation ``$\ll_M$" is always transitive.
For subsets $U,V\subset M$, we write $U\leq_M V$ ($U\ll_M V$ and $U<_M V$, respectively) if $x\leq_M y$ ($x\ll_M y$ and $x<_M y$, respectively) for any $x\in U$ and $y\in V$.
Here, we also show the {\it openness} of the relation ``$\ll_M$", i.e., for any two points $x,y\in M$, if $x\ll_M y$, there exist neighborhoods $U_1$ of $x$ and $V_1$ of $y$ such that $U_1\ll_M V_1$. See \cite[Proposition 2.8]{NW24}.

\begin{defn}\label{quasi-closed}
	The order ``$\leq_M$" is said to be \textit{quasi-closed} if $x\leq_M y$ whenever $x_n\to x$ and $y_n\to y$ as $n\to \infty$ and $x_n\ll_M y_n$ for all $n$.
\end{defn}

\begin{rmk}\label{rk:not-clo}
	In general, the conal order is not a closed order (the order ``$\leq_M$" is said to be \textit{closed} if $x\leq_M y$ whenever $x_n\to x$ and $y_n\to y$ as $n\to \infty$ and $x_n\leq_M y_n$), see \cite[p.299]{Lawson89} and \cite[p.470]{Neeb91}. The quasi-closed order relationship here is inspired by the causal continuity of space-times in \cite{Hawkingandsachs74}.
	Refer to \cite{NW24} for more details.
\end{rmk}

The conal order ``$\leq_M$" is called a \textit{partial} order relation if it is additionally antisymmetric (i.e., $x=y$ whenever $x\leq_M y$ with $y\leq_M x$).
A partial order ``$\leq$" is said to be \textit{locally order convex} if at each point, it has a basis of neighborhoods that are order convex ($z, x\in U$ implies $\{y: z\leq y \leq x\}\subset U$). See \cite{Lawson89,Neeb91}.

\begin{defn}\label{globally orderable}
    A conal manifold is \textit{globally orderable} if the order ``$\leq_M$" is a partial order which is locally order convex.
\end{defn}

\begin{rmk}
	The question of globality --- whether the conal order associated with a cone field extends to a partial order on the entire manifold --- is a central problem studied in \cite{Lawson89}. It is shown there that, on a homogeneous manifold, the globality of the conal order is equivalent to the globality of the Lie wedge in the corresponding Lie group action (see \cite[Section 5]{Lawson89}, \cite[Theorem 1.6]{Neeb91}, or \cite[Section 4.3]{HilgertNeeb93}). Global conal orders arise naturally in various contexts. For instance, the conal order induced by the cone field $C_{S^{+}_{n}}(A)$ on $S^{+}_{n}$, as mentioned above, defines a partial order (see \cite[Theorem 2]{MostajeranSepulchre18orderingmatrices}). Moreover, if a manifold is globally orderable, then closed conal curves can be excluded (see \cite[Section 5]{Lawson89}).
\end{rmk}

\vskip 4mm

\begin{defn}\label{measurezerodefi}
	If $M$ is a smooth manifold of dimension $r$, we say that a subset $A\subset M$ has \textit{measure zero} in $M$ if for every smooth chart $(U,\Phi)$ for $M$, the subset $\Phi(A\cap U)\subset \mathbb{R}^r$ has measure zero.
\end{defn}

\begin{rmk}\label{countabletest}
	This concept of sets of measure zero is independent of the metric of $M$.
	Moreover, to verify whether a set $A\subset M$ has measure zero, it suffices to check this property in a countable collection of smooth charts whose domains cover $A$ (see \cite[Lemma 6.6]{LeeGTM218}).
\end{rmk}


In the present paper, $M$ is a Riemannian manifold giving rise to a canonical measure structure (namely, the Riemann-Lebesgue measure). Recall that a subset $A$ of $M$ is said to be \textit{measurable} if around every $p\in A$ there is a chart $(U,\Phi)$ such that $\Phi(A\cap U)$ is (Lebesgue) measurable in $\mathbb{R}^r$. The \textit{Lebesgue $\sigma$-algebra of $M$} is denoted by $\mathcal{L}_M\triangleq \{A\subset M : A \text{ is measurable}\}$. Clearly, $\mathcal{L}_M$ contains the Borel $\sigma$-algebra of $M$.
Let $g$ be the Riemannian metric on $M$, and let $(U, \Phi)$ be a chart of $M$ with $\Phi=(x^1,\cdots,x^r)$, then the push forward $\Phi_{*}$ and the Gram determinant $G=\det [g_{jk}]$ with $g_{jk}=g(\frac{\partial}{\partial x^j},\frac{\partial}{\partial x^k})$ are well defined.
For $A\in \mathcal{L}_M$ with $A\subset U$, we define $\text{vol}_{g,U}(A)\triangleq \int _{\Phi(A)}\Phi_{*}\sqrt{G}dx$ (see \cite{AmannandEscher} for more details).

For $A\in \mathcal{L}_M$, the \textit{volume} of $A$ is denoted by
\[
\text{vol}_{g}(A)\triangleq \sum_{j=0}^{\infty}\text{vol}_{g,U_j}(A_j),
\]
where $\{(\Phi_{j}, U_j): j\in \mathbb{N}\}$ is a countable atlas of $M$, $A_0\triangleq A\cap U_0$ and $A_{n+1}\triangleq (A\cap U_{n+1})\backslash \cup_{k=0}^{n}A_{k}$ for $n\in \mathbb{N}$.

The \textit{Riemann-Lebesgue measure} of $M$ is denoted by $\lambda_M\triangleq \text{vol}_{g}$. We say that a set $A\subset M$ \textit{has full Riemann–Lebesgue measure} if the complement of $A$ in $M$ is a $\lambda_M$-null set.

\begin{rmk}\label{measueRMK}
	A measurable subset $A\subset M$ is a $\lambda_M$-null set, if and only if $A$ has measure zero in $M$ (see Proposition 1.6 in \cite[Chapter XII]{AmannandEscher}).
\end{rmk}

\vskip 4mm

Let $\Sigma$ be a dynamical system on $M$, generated by a continuously differentiable vector field $f$, where $M$ is equipped with a continuous solid cone field $C_M$. The flow of $\Sigma$ is denoted by $\varphi \colon \mathbb{R} \times M \to M$, $(t, x) \mapsto \varphi_t(x)$. For brevity, we also denote the flow by $\varphi_t$. Let $d\varphi_t(x)$ be the tangent map from $T_{x}M$ to $T_{\varphi_{t}(x)}M$.
For $x\in M$, the \textit{positive semiorbit} (\textit{negative semiorbit}, respectively) of $x$, be denoted by $O^+(x)=\{\varphi_t(x):t\geq 0\}$ ($O^-(x)=\{\varphi_t(x):t\leq 0\}$, respectively). The \textit{full orbit} of $x$ will be denoted by $O(x)=O^+(x)\cup O^-(x)$. An \textit{equilibrium} is a point $x$ for which $O(x)=\{x\}$. The set of equilibria is denoted by $\mathcal{E}$. The omega limit set $\omega(x)$ of $x$ is defined by $\omega(x)=\cap_{s\geq 0}\overline{\cup_{t\geq s}\varphi_t(x)}$. A point $y\in \omega(x)$ if and only if there is a sequence $\{t_i\}$, $t_i\to \infty$, such that $\varphi_{t_i}(x)\to y$ as $i\to \infty$. If $O^+(x)$ possesses a compact closure, then $\omega(x)$ is non-empty, compact, connected, and invariant ($\varphi_t(\omega(x))=\omega(x)$ for any $t\in  \mathbb{R}$). A point $x$ is a \textit{convergent point} if $\omega(x)$ consists of a single point of $\mathcal{E}$. The set of all convergent points is denoted by $\mathcal{C}$.

Throughout this paper, we always assume that each orbit of $\Sigma$ is well-defined for all $t\in \mathbb{R}$ and possesses a compact closure.


\vskip 4mm

\begin{defn}\label{uniformlySDP}
\begin{enumerate}[(i)]
		\item The system $\Sigma$ is \textit{differentially positive} (DP) with respect to $C_{M}$ if
        \begin{equation*}
        d\varphi_t(x)C_{M}(x)\subseteq C_{M}(\varphi_t(x)), \ \ \forall x\in M, \ \ \forall t\geq 0.
        \end{equation*}
        \item The differentially positive system $\Sigma$ is \textit{strongly differentially positive} (SDP) with respect to $C_{M}$ if
        \begin{equation*}
     	d\varphi_t(x)\{ C_{M}(x)\backslash \{ 0\} \} \subseteq \text{Int} C_{M}(\varphi_t(x)), \ \ \forall x\in M, \ \ \forall t>0.
        \end{equation*}
        \item The system $\Sigma$ is \textit{uniformly contracting} with respect to $C_{M}$ (which is $\Gamma$-invariant) if there exists $t_0>0$ and a cone field $R_{M}(x)\subseteq {\rm Int } C_M(x) \cup \{0\}$ on $M$ such that {\rm (a)} $d\varphi_t(x)C_{M}(x)\subseteq R_{M}(\varphi_t(x))$ for all $x\in M$ and all $t\geq t_0$; {\rm (b)} $\Gamma(x_1,x_2)R_M(x_1)=R_M(x_2)$ for each $x_1,x_2\in M$, where $\Gamma$ is defined above.
\end{enumerate}
\end{defn}

\begin{rmk}
	The uniformly contracting property in Definition \ref{uniformlySDP}(iii) implies a generalization of Birkhoff's contraction on $M$ (see Forni-Sepulchre \cite[Theorem 2]{ForniandSepulchre16} and its detailed proof in \cite[Appendix A]{ForniandSepulchre16}). Furthermore, if each orbit of $\Sigma$ is well-defined for all $t\in \mathbb{R}$, such uniform contraction ensures the existence of \textit{Perron-Frobenius vector field} (see \cite[Section 4]{NW24}, or \cite[Theorem 3]{ForniandSepulchre16} and \cite[Appendix A]{ForniandSepulchre16}).
\end{rmk}

In the present paper, we will prove Conjecture A for system $\Sigma$ which is (SDP) and uniformly contracting. Before ending this section, we impose the following hypotheses:
\begin{enumerate}[{\bf (H1)}]
	\item $M$ is a globally orderable conal manifold equipped with a continuous solid cone field $C_M$.
	\item The conal order ``$\leq_M$" is quasi-closed.
	\item Both the cone field $C_M$ and the Riemannian metric on $M$ are $\Gamma$-invariant.
\end{enumerate}

\section{Preliminary results}

Before presenting our main result, we give several preliminary results, which are crucial to our approach.

The following proposition indicates that differentially positive system $\Sigma$ naturally preserves the conal order.

\begin{prop}\label{strong order preserving}
	Assume that the system $\Sigma$ is {\rm (DP)} and uniformly contracting. Then, if $x<_M y$, one has
	\begin{enumerate}[{\rm (1)}]
		\item $\varphi_{t}(x)<_M \varphi_{t}(y)$ for all $t>0$;
		\item $\varphi_{t}(x)\ll_M \varphi_{t}(y)$ for all $t\geq t_0$;
		\item there exist neighborhoods $U$ of $x$, $V$ of $y$ such that $\varphi_{t_0}(U)\ll_M \varphi_{t_0}(V)$, where $t_0$ is defined in Definition {\rm \ref{uniformlySDP}(iii)}.
	\end{enumerate}
\end{prop}
\begin{proof}
	If $x<_M y$, there exists a conal curve $\gamma(s)$ such that $\gamma(0)=x$ and $\gamma(1)=y$. For (1), since $\Sigma$ is differentially positive and $\frac{d}{ds}\gamma(s)\in C_{M}(\gamma(s))\backslash \{ 0\}$, then $\frac{d}{ds}\varphi_{t}(\gamma(s))=d\varphi_{t}(\gamma(s))\frac{d}{ds}\gamma(s)\in C_{M}(\varphi_t(\gamma(s)))$ for all $t>0$. So, $\varphi_{t}(\gamma(s))$ is a conal curve and $\varphi_{t}(x)<_M \varphi_{t}(y)$ for all $t>0$. For {\rm (2)}, the fact that $d\varphi_t(x)C_{M}(x)\subseteq R_{M}(\varphi_t(x))\subseteq \text{Int} C_{M}(\varphi_{t}(\gamma(s)))$ for all $t\geq t_0$ implies that the conal curve $\varphi_{t}(\gamma(s))$ is a strictly conal curve and $\varphi_{t}(x)\ll_M \varphi_{t}(y)$ for all $t\geq t_0$, where $t_0$ is defined in Definition \ref{uniformlySDP}(iii). Finally, for {\rm (3)}, together with that $\varphi_{t_0}(x)\ll_M \varphi_{t_0}(y)$, the openness of the relation ``$\ll_M$" implies that there exist neighborhoods $\bar{U}$ of $\varphi_{t_0}(x)$ and $\bar{V}$ of $\varphi_{t_0}(y)$ such that $\bar{U}\ll_M\bar{V}$. By the continuity of $\varphi_{t_0}$, there are neighborhoods $U$ of $x$, $V$ of $y$ such that $\varphi_{t_0}(U)\subset \bar{U}$ and $\varphi_{t_0}(V)\subset \bar{V}$. Thus, we have proved the proposition.
\end{proof}
\vskip 2mm

We then give the following critical lemma, which turns out to be important for the proof of our main results in the forthcoming section.

\begin{lem}\label{fundamentalresults}
	Assume that {\rm (H1)-(H3)} hold. Assume also that $\Sigma$ is {\rm (SDP)} and uniformly contracting. Then, one has
	\begin{enumerate}[{\rm (a)}]
		\item{\rm (Non-ordering of limit sets).} The $\omega$-limit set cannot contain two points $x$ and $y$ with $x<_M y$;
		
		\item{\rm (Limit set dichotomy).} If $x<_M y$, then either $\omega(x)\ll_M \omega(y)$, or $\omega(x)=\omega(y)=\{e\}$ for some $e\in \mathcal{E}$.
	\end{enumerate}
\end{lem}

\begin{proof}
	This lemma is a combination of several results in \cite{NW24} (see \cite[Lemma 3.1, Lemma 4.6 and Proposition 4.7]{NW24}). More precisely, {\rm (a)} follows from \cite[Lemma 3.1(a)]{NW24} and our Proposition \ref{strong order preserving}. For {\rm (b)}, if $x<_M y$, it follows from \cite[Lemma 3.1(b)]{NW24} that
	\begin{equation}\label{dichotomy-equ}
		\text{either } \omega(x)\ll_M \omega(y) \text{ or } \omega(x)=\omega(y)\subset \mathcal{E}.
	\end{equation}
	Moreover, for the improved version as $\omega(x)=\omega(y)=\{e\}$ for some $e\in \mathcal{E}$, one can see \cite[Lemma 4.6 and Proposition 4.7]{NW24}.
\end{proof}

\begin{rmk}
	As we know, the conal order ``$\leq_M$" is in general not a closed relation (see Remark \ref{rk:not-clo} in Section $2$), due to the topological and geometric structure of the manifold $M$. However, as pointed by Hirsch \cite[p.28, line 9]{H88}, the closedness of the order is crucial for the dichotomy \eqref{dichotomy-equ} in the classical theory. Consequently, in order to overcome the difficulty for the absence of closedness, we introduce (SDP) which ensures that ``$\varphi_{t}(x)\ll_M \varphi_{t}(y)$ whenever $x<_M y$ and $t\in (0,t_0)$ for some $t_0>0$". Together with the uniformly contracting property, we can obtain the dichotomy \eqref{dichotomy-equ} (see more details in \cite[Lemma 3.6, Proposition 3.8 and the proof of Lemma 3.1(b)]{NW24}).
\end{rmk}

\begin{rmk}
	Of course, when the conal order ``$\leq_M$" is closed (see, e.g., \cite{MostajeranSepulchre18homogeneous,MostajeranSepulchre18orderingmatrices}), according to our approach in \cite{NW24}, Lemma \ref{fundamentalresults} remains valid provided that $\Sigma$ is (DP) and uniformly contracting. Such $\Sigma$ is also called as uniformly strictly differentially positive in Forni-Sepulchre \cite{ForniandSepulchre16}.
\end{rmk}

\begin{rmk}
	The uniformly contracting assumption in Lemma \ref{fundamentalresults} guarantees the existence of Perron-Frobenius vector field on $M$ (see \cite[Section 4]{NW24}, or \cite[Theorem 3]{ForniandSepulchre16} and \cite[Appendix A]{ForniandSepulchre16}), which turns out to be a key tool for the works in \cite{NW24,ForniandSepulchre16}.
\end{rmk}

\section{Main result on almost sure convergence}

In this section, we present our main result, which gives an affirmative answer to Conjecture A for the almost sure convergence of differentially positive system $\Sigma$. More precisely, we have

\begin{thm}[Almost Sure Convergence]\label{measure zero}
	Assume that {\rm (H1)-(H3)} hold. Assume also that system $\Sigma$ is {\rm (SDP)} and uniformly contracting. Then the set $\mathcal{C}$ of convergent points has full Riemann-Lebesgue measure.
\end{thm}

Before proceeding our proof, we give several technical lemmas, which turn out to be important to our approach.
Recall that $\mathcal{C}$ denotes the set of all convergent points.

\begin{lem}\label{borel set}
	The set $\mathcal{C}$ of convergent points is a Borel set.
\end{lem}
\begin{proof}
	Let $\{b_{i}\}_{i\in \mathbb{N}}$ be a dense countable collection of elements of $M$ and $$B_{k^{-1}}(b_{i})=\{x\in M : d(b_{i},x)\leq k^{-1}\}.$$
	Denote  $$W(B_{k^{-1}}(b_{i}),r)=\{x\in M : \varphi_{t}(x)\in B_{k^{-1}}(b_{i}) \text{ for all } t\geq r\}.$$
	Then $W(B_{k^{-1}}(b_{i}),r)$ is a Borel set since $W(B_{k^{-1}}(b_{i}),r)=\cap _{q\in R,\ q>r} \varphi_{q}^{-1}(B_{k^{-1}}(b_{i}))$, where $R$ denotes the set of rational numbers.
	
	Then the statement follows from the following equation:
	$$\mathcal{C}=\cap_{k\in \mathbb{N}}\cup_{i,r\in \mathbb{N}}W(B_{k^{-1}}(b_{i}),r).$$
	To prove it, first let $x\in \mathcal{C}$. Since for any fixed $k$ there exists $b_{i}$ such that $d(b_{i},\omega(x))\leq\frac{1}{2k}$ and therefor $x\in W(B_{k^{-1}}(b_{i}),r)$ for some large enough $r$, it has $x\in \cap_{k\in \mathbb{N}}\cup_{i,r\in \mathbb{N}}W(B_{k^{-1}}(b_{i}),r).$ Conversely, let $x\in \cap_{k\in \mathbb{N}}\cup_{i,r\in \mathbb{N}}W(B_{k^{-1}}(b_{i}),r).$ For each $k$, let $a_{k}$, $r_{k}$ be such that $x\in W(B_{k^{-1}}(a_{k}),r_{k})$. Then, $\{a_{k}\}$ is a Cauchy sequence; and therefore, converges to a point $a\in M$. Given $\varepsilon>0$, choose $k$ sufficiently large such that $d(\varphi_{r}(x),a)<\varepsilon$ for any $r\geq r_{k}$. Thus, we have $\omega(x)=a$ and therefor $x\in \mathcal{C}$.
\end{proof}

%

\begin{lem}\label{countable set}
	Assume that {\rm (H1)-(H3)} hold and system $\Sigma$ is {\rm (SDP)} and uniformly contracting.
	Let $L\subset M$ be a simple ordered set {\rm (}i.e., for any $p\neq q\in L$, either $p<_M q$ or $q<_M p${\rm )} of points. Then, $L\backslash \mathcal{C}$ is countable.
\end{lem}
\begin{proof}
We first claim that for any $x\in L\backslash \mathcal{C}$, there exists an open neighborhood $U(x)$ of $\omega(x)$ such that $U(x)\cap \omega(y)= \emptyset$ for every $y\in L\backslash \{x\}$.

Before proving the claim, we show how it implies this lemma. In fact, by virtue of the axiom of choice, one can define a mapping
	\[
	h : L\backslash \mathcal{C} \to M, \,\, x \mapsto h(x)\triangleq p_x,
	\]
by choosing some $p_x \in \omega(x)$ for each $x\in L\backslash \mathcal{C}$. Clearly, $h$ is injective (In fact, for any $x, z \in L\backslash \mathcal{C}$, since $p_x \in \omega(x)$ and $p_{z} \in \omega(z)$, the claim entails that $p_x\neq p_{z}$ whenever $x\neq z$). Now, choose a countable (neighborhood) base $\mathcal{B}$ of $M$ (with $\mathcal{B}\subset 2^M$). Again, by the axiom of choice, one can find a set-valued mapping $g$, associated with the above mapping $h$, as follows:
	\[
	g : L\backslash \mathcal{C} \to \mathcal{B}, \ \ \ \ x \mapsto g(x)\triangleq V(x),
	\]
with $V(x)\in \mathcal{B}$ satisfying $p_x\in V(x)\subset U(x)$. Here, the neighborhood $U(x)$ is as in the claim, and $p_x$ is defined in the map $h$. So, $g$ is also injective (In fact, the claim directly implies that $U(x)\cap \omega(z)= \emptyset$ whenever $x\neq z$. Hence, $V(x)\cap \omega(z)= \emptyset$. Therefore, by noticing that $p_{z}\in \omega(z)$, one has $p_{z}\notin V(x)$. This implies that $V(x)\neq V(z)$ whenever $x\neq z$, which means that $g$ is injective). Consequently, together with the fact that $\mathcal{B}$ is a countable, we then obtain that $L\backslash \mathcal{C}$ is countable, which completes the proof of the lemma.
	
Now, it remains to prove the claim. For $x\in L\backslash \mathcal{C}$, we define
$$Y_+(x)=\bigcup \{\omega(y) : x<_M y, \,y\in L\}.$$  Noticing that $x\notin \mathcal{C}$, it then follows from Lemma \ref{fundamentalresults}(b) (i.e., the limit set dichotomy) that $\omega(x)\ll_M \omega(y)$ for any $x<_M y$. Hence, $\omega(x) \ll_M Y_+(x)$, for any $x\in L\backslash \mathcal{C}$. So, by the quasi-closedness of the conal order ``$\leq_M$" in (H2), we obtain $\omega(x) \leq_M \overline{Y_+(x)}$. This implies that $\overline{Y_+(x)}\cap \omega(x)= \emptyset$ (For otherwise, we can choose $p\in \overline{Y_+(x)}\cap \omega(x)\ne \emptyset$. So,
$\omega(x) \leq_M p$. Since $x\notin \mathcal{C}$, there exists some $q\in \omega(x)$ such that $q<_M p$, which contradicts  Lemma \ref{fundamentalresults}(a) (i.e., non-ordering of the limit sets)).

Similarly, for $x\in L\backslash \mathcal{C}$, we can define $$Y_-(x)=\bigcup \{\omega(y) : y<_M x, \,y\in L\},$$ and obtain that $\overline{Y_-(x)}\cap \omega(x)= \emptyset$. Consequently, we have
	\[
	(\overline{Y_-(x)\cup Y_+(x)})\cap \omega(x)=\emptyset,\,\, \text{ for any } x\in L\backslash \mathcal{C}.
	\]
	Since $M$ is a complete metric space, there exist open neighborhoods $U(x)$ and $W(x)$ of $\omega(x)$ and $\overline{Y_-(x)\cup Y_+(x)}$, respectively, such that $U(x)\cap W(x)= \emptyset$. Noticing that $\omega(y)\subset \overline{Y_-(x)\cup Y_+(x)}$ for any $y\in L\backslash \{x\}$, then it directly yields that $U(x)\cap \omega(y)= \emptyset$. Thus, we have proved the claim. This completes the proof.
\end{proof}

Now, we are ready to prove Theorem \ref{measure zero}.

\begin{proof}[Proof of Theorem \ref{measure zero}]
Let $D\triangleq M\backslash \mathcal{C}$. By virtue of Remark \ref{measueRMK} with Lemma \ref{borel set}, one only needs to prove that $D$ has measure zero in $M$.

Fix an $x\in M$ and suppose $(U,\Phi)$ be a smooth chart around $x$. Let $\tilde{v}\in \text{Int}C_{M}^{\Phi}(x)$ be a unit vector in $\mathbb{R}^{r}$.
We decompose $\mathbb{R}^{r}$ directly into $E$ and $F$, where $E=\{l\tilde{v}: l\in \mathbb{R}\}$ and $F\subset\mathbb{R}^{r}$ is a closed linear subspace complementary to $E$, i.e., $\mathbb{R}^{r}=E\oplus F$. Let $\mu_E$ and $\mu_F$ denote the Lebesgue measures on $E$ and $F$, respectively. Then one can define the measure $\mu$ on $\mathbb{R}^{r}$ as $\mu(A\times B)=\mu_{E}(A)\times\mu_{F}(B)$ for any set $A\subset E$ and $B\subset F$ (whose minimal $\sigma$-algebra is the all Borel sets in $\mathbb{R}^{r}$). Clearly, $\mu$ is the Lebesgue measure on $\mathbb{R}^{r}$ (see, e.g., \cite[Chapter II]{HalmosGTM18}). For each $\tilde{y}\in F$, we write $E_{\tilde{y}}=\tilde{y}+E$, and let $\mu_{\tilde{y}}$ be the Lebesgue measure on the $E_{\tilde{y}}$. By Fubini's Theorem, it has
		\begin{equation}\label{fubini}
		\mu(A)=\int_{\tilde{y}\in F}\mu_{\tilde{y}}(E_{\tilde{y}}\cap A)d \mu_{F},\, \ \text{for any Borel set } A\subset \mathbb{R}^{r}.
	\end{equation}
Since the cone field is continuous, one can choose a small neighborhood $V(\subset U)$ of $x$ such that $\tilde{v}\in \text{Int}C_{M}^{\Phi}(z)$ for any $z\in V$. Let $\tilde{V}\triangleq \Phi (V) \subseteq\mathbb{R}^{r}$ and $\tilde{D}\triangleq\Phi (D\cap U) \subseteq\mathbb{R}^{r}$. It then follows from Lemma \ref{borel set} that $\tilde{D}$ is a Borel set.

We now {\it claim that $\Phi^{-1}(E_{\tilde{y}}\cap \tilde{V})$ is ordered for any $\tilde{y}\in \Pi_{F}(\tilde{V})$}, where $\Pi_{F}: \mathbb{R}^{r}\rightarrow F$ be the projection with kernel $E$. In fact, for any fixed $\tilde{y}\in \Pi_{F}(\tilde{V})$, $E_{\tilde{y}}\cap \tilde{V}$ is a line segment in $\mathbb{R}^{r}$ whose tangent vector at each point has the same direction as the vector $\tilde{v}$. Since $\tilde{v}\in \text{Int}C_{M}^{\Phi}(z)$ for any $z\in \Phi^{-1}(E_{\tilde{y}}\cap \tilde{V})\subset V$, one obtains that the tangent vector at each point $\tilde{p}$ of $E_{\tilde{y}}\cap \tilde{V}$ lies in the cone $C_{M}^{\Phi}(\Phi^{-1}(\tilde{p}))$. Thus, $\Phi^{-1}(E_{\tilde{y}}\cap \tilde{V})$ is ordered, which completes the proof of the claim.
	
Due to the claim, by taking $L$ as $\Phi^{-1}(E_{\tilde{y}}\cap \tilde{V})$ in Lemma \ref{countable set}, we obtain that for each $\tilde{y}\in \Pi_{F}(\tilde{V})$, the set $\Phi^{-1}(E_{\tilde{y}}\cap \tilde{V})\cap D$ is at most countable in $M$. Hence, $E_{\tilde{y}}\cap \tilde{V}\cap \tilde{D} $ is at most countable in $\mathbb{R}^{r}$. This implies that $\mu_{\tilde{y}}(E_{\tilde{y}}\cap \tilde{V}\cap \tilde{D})=0$, for any $\tilde{y}\in \Pi_{F}(\Phi(V))$. Consequently, the integrand in (\ref{fubini}) is zero, i.e., $\mu(\tilde{V}\cap \tilde{D})=0$. Recall that $\tilde{V}\cap \tilde{D}=\Phi(V\cap D)$. Then, it has $\mu(\Phi(D\cap V))=0$. Here, it deserves to point that $V$ and $\Phi$ are dependent on $x$, and thus labeled by $V_{x}$ and $\Phi_{x}$, respectively. Accordingly, it induces that $\mu(\Phi_{x}(D\cap V_{x}))=0$ for any $x\in M$.
	
By the arbitrariness of $x\in M$, we have obtained a collection of open sets $\{V_{x}\}_{x\in M}$ cover $M$ satisfying $\mu(\Phi_{x}(D\cap V_{x}))=0$ in $\mathbb{R}^{r}$, for any $x\in M$. Since $M$ is second countable, there exists a countable open sets $\{V_{i}\}_{i\in \mathbb{N}}\subset \{V_{x}\}_{x\in M}$ cover $M$ (see \cite[Proposition A.16]{LeeGTM218}), which forms smooth charts $(V_{i},\Phi_{i})$ whose domains cover $M$ (see \cite[Lemma 1.35]{LeeGTM218}). Consequently, we have constructed a smooth charts $(V_{i},\Phi_{i})$ of $M$ such that $\mu(\Phi_{i}(D\cap V_{i}))=0$ in $\mathbb{R}^{r}$, for $i=1,2,\cdots$. Thus, by Remark \ref{countabletest}, $D$ has measure zero in $M$. We have completed the proof.
\end{proof}


\begin{rmk}\label{rkofthm}
Our Theorem \ref{measure zero}, in conjunction with \cite[Theorem 4.1]{NW24}, demonstrates that for differentially positive systems, {\it the set $\mathcal{C}$ is not only ``generic" in the topological sense, but also ``almost everywhere"  in the measure-theoretic sense}.
We further point out that the common idea is to figure out the characteristics of the points that are not convergent. For {\it ``the genericity"}, a key insight is to show that, if $x\notin \mathcal{C}$ then it must belong to the closure of $\mathcal{C}$ (i.e., $x\in \overline{\mathcal{C}}$, see \cite[Theorem 4.1]{NW24}). In other words, any such $x$ can be approached topologically by the points in $\mathcal{C}$. However, we stress that such observation is not enough to guarantee the ``almost sure convergence" at all.
In order to obtain our Theorem \ref{measure zero}, a more critical measure-theoretic insight is needed: in a simply ordered set, at most countably many points can be non-convergent (Clearly, this implies that any such point $x$ automatically belongs to $\overline{\mathcal{C}}$). Together with Fubini's Theorem, we succeed in establishing that the set of non-convergent points in each local chart has measure zero; and hence, the set of non-convergent points has measure zero in $M$.
\end{rmk}

\section*{Acknowledgments}
The authors are greatly indebted to Professor Karl-Hermann Neeb and Professor Wenxian Shen for very valuable discussions and fruitful explanations about conal orders for our present series of papers.

\end{document}